\documentclass[12pt]{article}
\usepackage{amsmath,amssymb,amsthm}
\usepackage{hyperref}
\usepackage{datetime}
\usepackage{units}
\usepackage{color}
\usepackage[T1]{fontenc}
\usepackage[utf8]{inputenc}
\usepackage{authblk}
\usepackage{bm,latexsym,mathrsfs,enumerate}
\usepackage{color}

\title{ A class of digit extraction BBP-type formulas in general binary bases\thanks{%
MSC 2010: 11Y60, 30B99}}
\author[1]{Kunle Adegoke\thanks{adegoke00@gmail.com}}
\author[2]{Jaume Oliver Lafont}
\author[3]{Olawanle Layeni}
\affil{Department of Physics and Engineering Physics, \mbox{Obafemi Awolowo University, Ile-Ife, 220005 Nigeria}}
\affil[2]{Conselleria d'Educacio i Cultura, Govern de les Illes Balears, Palma, Spain}
\affil[3]{Department of Mathematics, \mbox{Obafemi Awolowo University, Ile-Ife, 220005 Nigeria}}

\theoremstyle{plain}
\numberwithin{equation}{section}

\begin{document}
\date{}
\maketitle
\begin{abstract}
\noindent BBP-type formulas are usually discovered experimentally, one at a time and in specific bases, through computer searches. In this paper, however, we derive directly, without doing any searches, explicit digit extraction BBP-type formulas in general binary bases $b=2^{12p}$, for $p$ positive odd integers. As particular examples, new binary formulas are presented for $\pi\sqrt 3$, $\pi\sqrt 3\log 2$, $\sqrt 3\;{\rm Cl}_2(\pi/3)$ and a couple of other polylogarithm constants. A variant of the formula for $\pi\sqrt 3\log 2$ derived in this paper has been known for over ten years but was hitherto unproved.  Binary BBP-type formulas for the logarithms of an infinite set of primes and binary BBP-type representations for the arctangents of an infinite set of rational numbers are also presented. Finally, new binary BBP-type zero relations are established.
\end{abstract}
\tableofcontents

\section{Introduction}

BBP-type formulas are formulas of the form
\[
\alpha = \sum_{k=0}^\infty  1/b^k \sum_{j=1}^l a_j / (k l + j)^s
\]
where $s$, $b$, $l$ (degree, base, length respectively) and $a_j$ are integers, and $\alpha$ is some constant. Formulas of this type were first introduced in a 1996 paper~\cite{bbp97}, where a formula of this type for $\pi$ was given. Such formulas allow digit extraction --- the $i$-th digit of a mathematical constant $\alpha$ in base $b$ can be calculated directly, without needing to compute any of the previous $i-1$ digits, by means of simple algorithms that do not require multiple-precision arithmetic~\cite{bailey09}.

\bigskip

Apart from digit extraction, another reason the study of BBP-type formulas has continued to attract attention is that BBP-type constants are conjectured to be either rational or normal to base $b$~\cite{bailey01,borwein02,chamberland03}, that is their base-$b$ digits are randomly distributed. 

\bigskip

BBP-type formulas are usually discovered experimentally, one at a time and in specific bases, through computer searches. In this paper, however, we derive explicit digit extraction BBP-type formulas in general binary bases $b=2^{12p}$, for $p\in\mathbb{Z^+}$ and$\mod(p,2)=1$.

\section{Definitions and Notation}

The polylogarithm functions denoted by {\rm Li} in this paper are defined by
\[
{\rm Li}_s [z] = \sum\limits_{k = 1}^\infty  {\frac{{z^k }}{{k^s }}},\quad |z|\le 1\,,s\in\mathbb{Z^+}\,. 
\]

For $|z|=1$ and $x\in[0,2\pi]$ we have
\[
\begin{split}
{\rm Li}_{2n} [e^{ix} ]& = {\rm Gl}_{2n} (x) + i{\rm Cl}_{2n} (x)\\
{\rm Li}_{2n + 1} [e^{ix} ]& = {\rm Cl}_{2n + 1} (x) + i{\rm Gl}_{2n + 1} (x)\,,
\end{split}
\]
where ${\rm Gl}$ and ${\rm Cl}$ are Clausen sums~\cite{lewin81} defined, for $n\in\mathbb{Z^+}$ by

\[
\begin{split}
{\rm Cl}_{2n} (x)& = \sum\limits_{k = 1}^\infty  {\frac{{\sin kx}}{{k^{2n} }}},\quad {\rm Cl}_{2n + 1} (x) = \sum\limits_{k = 1}^\infty  {\frac{{\cos kx}}{{k^{2n + 1} }}}\\
{\rm Gl}_{2n} (x) &= \sum\limits_{k = 1}^\infty  {\frac{{\cos kx}}{{k^{2n} }}},\quad {\rm Gl}_{2n + 1} (x) = \sum\limits_{k = 1}^\infty  {\frac{{\sin kx}}{{k^{2n + 1} }}}\,.
\end{split}
\]

We shall find the following formulas useful:

\[
\begin{split}
{\rm Gl}_{2n} (x) &= ( - 1)^{1 + [n/2]} 2^{n - 1} \pi ^n {\rm B}_n (x/2\pi )/n!\\
\frac{1}{{m^{n - 1} }}{\rm Cl}_n (mx) &= \sum\limits_{r = 0}^{m - 1} {{\rm Cl}_{n} (x + 2\pi r/m)}\,.
\end{split}
\]
Here $[n/2]$ denotes the integer part of $n/2$ and $ {\rm B}_n$ are the Bernoulli polynomials defined by

\[
{\frac {t{{\rm e}^{xt}}}{{{\rm e}^{t}}-1}}=\sum _{n=0}^{\infty }{
\frac {{\rm B}_n (x) {t}^{n}}{n!}}\,.
\]

In order to save space, we will give the BBP-type formulas using the compact P-notation~\cite{bailey09}:

\begin{equation}\label{equ.b98zkzg}
P(s,b,l,A) \equiv \sum\limits_{k = 0}^\infty  {\frac{1}{{b^k }}\sum\limits_{j = 1}^l {\frac{{a_j }}{{(kl + j)^s }}} }\,,
\end{equation}

where $s$, $b$ and $l$ are integers, and \mbox{$A = (a_1, a_2,\ldots, a_l)$} is a vector of integers.

\section{Degree $s$ Formulas}	

Using the identities

\begin{equation}
{\rm Re\,Li}_s \left[ {\frac{1}{{\sqrt 2 ^p }}\exp ix} \right] = \sum\limits_{k = 1}^\infty  {\left[ {\frac{1}{{\sqrt 2 ^{pk} }}\frac{1}{{k^s }}\cos kx} \right]} 
\end{equation}

and

\begin{equation}
{\rm Im\,Li}_s \left[ {\frac{1}{{\sqrt 2 ^p }}\exp ix} \right] = \sum\limits_{k = 1}^\infty  {\left[ {\frac{1}{{\sqrt 2 ^{pk} }}\frac{1}{{k^s }}\sin kx} \right]} 
\end{equation}

for $x\in \{\pi/12,5\pi/12,7\pi/12,11\pi/12\} $ and the fact that 
\[
\cos \left( {\frac{\pi }{{12}}} \right) = \frac{{\sqrt 3  + 1}}{{2\sqrt 2 }} = \sin \left( {\frac{{7\pi }}{{12}}} \right)= \sin \left( {\frac{{5\pi }}{{12}}} \right)=- \cos \left( {\frac{{11\pi }}{{12}}} \right)
\]
and
\[
\sin \left( {\frac{\pi }{{12}}} \right) = \frac{{\sqrt 3  - 1}}{{2\sqrt 2 }}=\cos \left( {\frac{{5\pi }}{{12}}} \right) =\sin \left( {\frac{{11\pi }}{{12}}} \right)=  - \cos \left( {\frac{{7\pi }}{{12}}} \right)\,,
\]
it is not difficult to obtain the following results, written in the P-notation (\eqref{equ.b98zkzg}):
\begin{eqnarray}
&& {\rm Re\,Li}_s \left[ {\frac{1}{{\sqrt 2 ^p }}\exp \left( {\frac{{i\pi }}{{12}}} \right)} \right] + {\rm Re\,Li}_s \left[ {\frac{1}{{\sqrt 2 ^p }}\exp \left( {\frac{{7i\pi }}{{12}}} \right)} \right] \nonumber\\ 
&& = \frac{1}{{2^{12p} }}P(s,2^{12p} ,24,(2^{( - {\textstyle{1 \over 2}} + {\textstyle{p \over 2}} + 11p)} ,0,2^{({\textstyle{1 \over 2}} - {\textstyle{p \over 2}} + 11p)} ,2^{10p} , \nonumber\\ 
&& - 2^{( - {\textstyle{1 \over 2}} + {\textstyle{p \over 2}} + 9p)} ,0,2^{( - {\textstyle{1 \over 2}} + {\textstyle{p \over 2}} + 8p)} , - 2^{8p} , - 2^{({\textstyle{1 \over 2}} - {\textstyle{p \over 2}} + 8p)} ,0, - 2^{( - {\textstyle{1 \over 2}} + {\textstyle{p \over 2}} + 6p)} , \nonumber\\ 
&&  - 2^{1 + 6p} , - 2^{( - {\textstyle{1 \over 2}} + {\textstyle{p \over 2}} + 5p)} ,0, - 2^{({\textstyle{1 \over 2}} - {\textstyle{p \over 2}} + 5p)} , - 2^{4p} ,2^{( - {\textstyle{1 \over 2}} + {\textstyle{p \over 2}} + 3p)} ,0, \nonumber\\ 
&&  - 2^{( - {\textstyle{1 \over 2}} + {\textstyle{p \over 2}} + 2p)} ,2^{2p} ,2^{({\textstyle{1 \over 2}} - {\textstyle{p \over 2}} + 2p)} ,0,2^{( - {\textstyle{1 \over 2}} + {\textstyle{p \over 2}})} ,2))\,,
 \end{eqnarray}
\begin{eqnarray}
&&{\rm Re\,Li}_s \left[ {\frac{1}{{\sqrt 2 ^p }}\exp \left( {\frac{{i\pi }}{{12}}} \right)} \right] - {\rm Re\,Li}_s \left[ {\frac{1}{{\sqrt 2 ^p }}\exp \left( {\frac{{7i\pi }}{{12}}} \right)} \right] \nonumber\\ 
&&\qquad  = \frac{{\sqrt 3 }}{{2^{12p} }}P(s,2^{12p} ,24,(2^{( - {\textstyle{1 \over 2}} + {\textstyle{p \over 2}} + 11p)} ,2^{11p} ,0,0, \nonumber\\ 
&& \qquad 2^{( - {\textstyle{1 \over 2}} + {\textstyle{p \over 2}} + 9p)} ,0, - 2^{( - {\textstyle{1 \over 2}} + {\textstyle{p \over 2}} + 8p)} ,0,0, - 2^{7p} , - 2^{( - {\textstyle{1 \over 2}} + {\textstyle{p \over 2}} + 6p)} , \nonumber\\ 
&&\qquad 0, - 2^{( - {\textstyle{1 \over 2}} + {\textstyle{p \over 2}} + 5p)} , - 2^{5p} ,0,0, - 2^{( - {\textstyle{1 \over 2}} + {\textstyle{p \over 2}} + 3p)} ,0, \nonumber\\ 
&&\qquad 2^{( - {\textstyle{1 \over 2}} + {\textstyle{p \over 2}} + 2p)} ,0,0,2^p ,2^{( - {\textstyle{1 \over 2}} + {\textstyle{p \over 2}})} ,0))\,,
\end{eqnarray}
\begin{eqnarray}\label{equ.aaob80z}
&& {\rm Im\,Li}_s \left[ {\frac{1}{{\sqrt 2 ^p }}\exp \left( {\frac{{i\pi }}{{12}}} \right)} \right] + {\rm Im\,Li}_s \left[ {\frac{1}{{\sqrt 2 ^p }}\exp \left( {\frac{{7i\pi }}{{12}}} \right)} \right] \nonumber\\ 
&&  = \frac{{\sqrt 3 }}{{2^{12p} }}P(s,2^{12p} ,24,(2^{( - {\textstyle{1 \over 2}} + {\textstyle{p \over 2}} + 11p)} ,0,0,2^{10p} , \nonumber\\ 
&& 2^{( - {\textstyle{1 \over 2}} + {\textstyle{p \over 2}} + 9p)} ,0,2^{( - {\textstyle{1 \over 2}} + {\textstyle{p \over 2}} + 8p)} ,2^{8p} ,0,0,2^{( - {\textstyle{1 \over 2}} + {\textstyle{p \over 2}} + 6p)} , \nonumber\\ 
&& 0, - 2^{( - {\textstyle{1 \over 2}} + {\textstyle{p \over 2}} + 5p)} ,0,0, - 2^{4p} , - 2^{( - {\textstyle{1 \over 2}} + {\textstyle{p \over 2}} + 3p)} ,0, \nonumber\\ 
&&  - 2^{( - {\textstyle{1 \over 2}} + {\textstyle{p \over 2}} + 2p)} , - 2^{2p} ,0,0, - 2^{( - {\textstyle{1 \over 2}} + {\textstyle{p \over 2}})} ,0))\,,
 \end{eqnarray}
\begin{eqnarray}\label{equ.escpvp9}
&& {\rm Im\,Li}_s \left[ {\frac{1}{{\sqrt 2 ^p }}\exp \left( {\frac{{i\pi }}{{12}}} \right)} \right] - {\rm Im\,Li}_s \left[ {\frac{1}{{\sqrt 2 ^p }}\exp \left( {\frac{{7i\pi }}{{12}}} \right)} \right] \nonumber\\ 
&&  = \frac{1}{{2^{12p} }}P(s,2^{12p} ,24,( - 2^{( - {\textstyle{1 \over 2}} + {\textstyle{p \over 2}} + 11p)} ,2^{11p} ,2^{({\textstyle{1 \over 2}} - {\textstyle{p \over 2}} + 11p)} ,0, \nonumber\\ 
&& 2^{( - {\textstyle{1 \over 2}} + {\textstyle{p \over 2}} + 9p)} ,2^{1 + 9p} ,2^{( - {\textstyle{1 \over 2}} + {\textstyle{p \over 2}} + 8p)} ,0,2^{({\textstyle{1 \over 2}} - {\textstyle{p \over 2}} + 8p)} ,2^{7p} , - 2^{( - {\textstyle{1 \over 2}} + {\textstyle{p \over 2}} + 6p)} , \nonumber\\ 
&& 0,2^{( - {\textstyle{1 \over 2}} + {\textstyle{p \over 2}} + 5p)} , - 2^{5p} , - 2^{({\textstyle{1 \over 2}} - {\textstyle{p \over 2}} + 5p)} ,0, - 2^{( - {\textstyle{1 \over 2}} + {\textstyle{p \over 2}} + 3p)} , - 2^{1 + 3p} , \nonumber\\ 
&&  - 2^{( - {\textstyle{1 \over 2}} + {\textstyle{p \over 2}} + 2p)} ,0, - 2^{({\textstyle{1 \over 2}} - {\textstyle{p \over 2}} + 2p)} , - 2^p ,2^{( - {\textstyle{1 \over 2}} + {\textstyle{p \over 2}})} ,0))\,, 
 \end{eqnarray}
\begin{eqnarray}
&& {\rm Re\,Li}_{\rm s} \left[ {\frac{1}{{\sqrt 2 ^p }}\exp \left( {\frac{{5i\pi }}{{12}}} \right)} \right] + {\rm Re\,Li}_{\rm s} \left[ {\frac{1}{{\sqrt 2 ^p }}\exp \left( {\frac{{11i\pi }}{{12}}} \right)} \right] \nonumber\\ 
&&= \frac{1}{{2^{12p} }}P\left( {s,2^{12p} ,24,\left( { - 2^{( - {\textstyle{1 \over 2}} + {\textstyle{p \over 2}} + 11p)} ,0,} \right.} \right. - 2^{({\textstyle{1 \over 2}} - {\textstyle{p \over 2}} + 11p)} ,2^{10p} , \nonumber\\ 
&& 2^{( - {\textstyle{1 \over 2}} + {\textstyle{p \over 2}} + 9p)} ,0, - 2^{( - {\textstyle{1 \over 2}} + {\textstyle{p \over 2}} + 8p)} , - 2^{8p} ,2^{({\textstyle{1 \over 2}} - {\textstyle{p \over 2}} + 8p)} ,0,2^{( - {\textstyle{1 \over 2}} + {\textstyle{p \over 2}} + 6p)} , \nonumber\\ 
&&  - 2^{1 + 6p} ,2^{( - {\textstyle{1 \over 2}} + {\textstyle{p \over 2}} + 5p)} ,0,2^{({\textstyle{1 \over 2}} - {\textstyle{p \over 2}} + 5p)} , - 2^{4p} , - 2^{( - {\textstyle{1 \over 2}} + {\textstyle{p \over 2}} + 3p)} ,0, \nonumber\\ 
&& \left. {\left. {2^{( - {\textstyle{1 \over 2}} + {\textstyle{p \over 2}} + 2p)} ,2^{2p} , - 2^{({\textstyle{1 \over 2}} - {\textstyle{p \over 2}} + 2p)} ,0, - 2^{( - {\textstyle{1 \over 2}} + {\textstyle{p \over 2}})} ,2} \right)} \right)\,, 
 \end{eqnarray}
\begin{eqnarray}
&& {\rm Re\,Li}_{\rm s} \left[ {\frac{1}{{\sqrt 2 ^p }}\exp \left( {\frac{{5i\pi }}{{12}}} \right)} \right] - {\rm Re\,Li}_{\rm s} \left[ {\frac{1}{{\sqrt 2 ^p }}\exp \left( {\frac{{11i\pi }}{{12}}} \right)} \right] \nonumber\\ 
&&  = \frac{{\sqrt 3 }}{{2^{12p} }}P\left( {s,2^{12p} ,24,\left( {2^{( - {\textstyle{1 \over 2}} + {\textstyle{p \over 2}} + 11p)} , - 2^{11p} ,} \right.} \right.0,0, \nonumber\\ 
&& 2^{( - {\textstyle{1 \over 2}} + {\textstyle{p \over 2}} + 9p)} ,0, - 2^{( - {\textstyle{1 \over 2}} + {\textstyle{p \over 2}} + 8p)} ,0,0,2^{7p} , - 2^{( - {\textstyle{1 \over 2}} + {\textstyle{p \over 2}} + 6p)} , \nonumber\\ 
&& 0, - 2^{( - {\textstyle{1 \over 2}} + {\textstyle{p \over 2}} + 5p)} ,2^{5p} ,0,0, - 2^{( - {\textstyle{1 \over 2}} + {\textstyle{p \over 2}} + 3p)} ,0, \nonumber\\ 
&& \left. {\left. {2^{( - {\textstyle{1 \over 2}} + {\textstyle{p \over 2}} + 2p)} ,0,0, - 2^p ,2^{( - {\textstyle{1 \over 2}} + {\textstyle{p \over 2}})} ,0} \right)} \right)\,, 
\end{eqnarray}
\begin{eqnarray}\label{equ.digtrb1}
&&{\rm Im\,Li}_{\rm s} \left[ {\frac{1}{{\sqrt 2 ^p }}\exp \left( {\frac{{5i\pi }}{{12}}} \right)} \right] + {\rm Im\,Li}_{\rm s} \left[ {\frac{1}{{\sqrt 2 ^p }}\exp \left( {\frac{{11i\pi }}{{12}}} \right)} \right] \nonumber\\ 
&& = \frac{{\sqrt 3 }}{{2^{12p} }}P\left( {s,2^{12p} ,24,\left( {2^{( - {\textstyle{1 \over 2}} + {\textstyle{p \over 2}} + 11p)} ,0,} \right.} \right.0, - 2^{10p} , \nonumber\\ 
&& 2^{( - {\textstyle{1 \over 2}} + {\textstyle{p \over 2}} + 9p)} ,0,2^{( - {\textstyle{1 \over 2}} + {\textstyle{p \over 2}} + 8p)} , - 2^{8p} ,0,0,2^{( - {\textstyle{1 \over 2}} + {\textstyle{p \over 2}} + 6p)} , \nonumber\\ 
&& 0, - 2^{( - {\textstyle{1 \over 2}} + {\textstyle{p \over 2}} + 5p)} ,0,0,2^{4p} , - 2^{( - {\textstyle{1 \over 2}} + {\textstyle{p \over 2}} + 3p)} ,0, \nonumber\\ 
&& \left. {\left. { - 2^{( - {\textstyle{1 \over 2}} + {\textstyle{p \over 2}} + 2p)} ,2^{2p} ,0,0, - 2^{( - {\textstyle{1 \over 2}} + {\textstyle{p \over 2}})} ,0} \right)} \right) 
\end{eqnarray}
and
\begin{eqnarray}\label{equ.uklodgq}
&&{\rm Im\,Li}_{\rm s} \left[ {\frac{1}{{\sqrt 2 ^p }}\exp \left( {\frac{{5i\pi }}{{12}}} \right)} \right] - {\rm Im\,Li}_{\rm s} \left[ {\frac{1}{{\sqrt 2 ^p }}\exp \left( {\frac{{11i\pi }}{{12}}} \right)} \right] \nonumber\\ 
&&= \frac{1}{{2^{12p} }}P\left( {s,2^{12p} ,24,\left( {2^{( - {\textstyle{1 \over 2}} + {\textstyle{p \over 2}} + 11p)} ,2^{11p} ,} \right.} \right. - 2^{({\textstyle{1 \over 2}} - {\textstyle{p \over 2}} + 11p)} ,0, \nonumber\\ 
&& - 2^{( - {\textstyle{1 \over 2}} + {\textstyle{p \over 2}} + 9p)} ,2^{1 + 9p} , - 2^{( - {\textstyle{1 \over 2}} + {\textstyle{p \over 2}} + 8p)} ,0, - 2^{({\textstyle{1 \over 2}} - {\textstyle{p \over 2}} + 8p)} ,2^{7p} ,2^{( - {\textstyle{1 \over 2}} + {\textstyle{p \over 2}} + 6p)} , \nonumber\\ 
&& 0, - 2^{( - {\textstyle{1 \over 2}} + {\textstyle{p \over 2}} + 5p)} , - 2^{5p} ,2^{({\textstyle{1 \over 2}} - {\textstyle{p \over 2}} + 5p)} ,0,2^{( - {\textstyle{1 \over 2}} + {\textstyle{p \over 2}} + 3p)} , - 2^{1 + 3p} , \nonumber\\ 
&& \left. {\left. {2^{( - {\textstyle{1 \over 2}} + {\textstyle{p \over 2}} + 2p)} ,0,2^{({\textstyle{1 \over 2}} - {\textstyle{p \over 2}} + 2p)} , - 2^p , - 2^{( - {\textstyle{1 \over 2}} + {\textstyle{p \over 2}})} ,0} \right)} \right)\,. 
 \end{eqnarray}

We note that although the above formulas hold true for all $p>0$, they are BBP-type only for $p\in\mathbb{Z^+}$ and$\mod(p,2)=1$.

\section{Degree $1$ Formulas}

When $s=1$, the polylogarithm constants on the left hand side in each of the above formulas can be evaluated, using the identities

\[
\arctan \left( {\frac{{q\sin x}}{{1 - q\cos x}}} \right) = {\rm Im\;Li}_1 \left[ {q\exp (ix)} \right]
\]

and

\[
 - \frac{1}{2}\log \left( {1 - 2q\cos x + q^2 } \right) = {\rm Re\;Li}_1 \left[ {q\exp (ix)} \right]\,,
\]
and we have the following degree $1$ binary BBP-type formulas
\begin{eqnarray}\label{equ.x9u9heu}
&& - \frac{1}{2}\log( {1 - 2^{^{({\textstyle{1 \over 2}} - {\textstyle{p \over 2}})} }  + 2^{ - p}  - 2^{^{({\textstyle{1 \over 2}} + {\textstyle{p \over 2}} - 2p)} }  + 2^{ - 2p} } ) \nonumber\\ 
&& \qquad  = \frac{1}{{2^{12p} }}P(1,2^{12p} ,24,(2^{( - {\textstyle{1 \over 2}} + {\textstyle{p \over 2}} + 11p)} ,0,2^{({\textstyle{1 \over 2}} - {\textstyle{p \over 2}} + 11p)} ,2^{10p} , \nonumber\\ 
&& \qquad  - 2^{( - {\textstyle{1 \over 2}} + {\textstyle{p \over 2}} + 9p)} ,0,2^{( - {\textstyle{1 \over 2}} + {\textstyle{p \over 2}} + 8p)} , - 2^{8p} , - 2^{({\textstyle{1 \over 2}} - {\textstyle{p \over 2}} + 8p)} ,0, - 2^{( - {\textstyle{1 \over 2}} + {\textstyle{p \over 2}} + 6p)} , \nonumber\\ 
&& \qquad  - 2^{1 + 6p} , - 2^{( - {\textstyle{1 \over 2}} + {\textstyle{p \over 2}} + 5p)} ,0, - 2^{({\textstyle{1 \over 2}} - {\textstyle{p \over 2}} + 5p)} , - 2^{4p} ,2^{( - {\textstyle{1 \over 2}} + {\textstyle{p \over 2}} + 3p)} ,0, \nonumber\\ 
&& \qquad  - 2^{( - {\textstyle{1 \over 2}} + {\textstyle{p \over 2}} + 2p)} ,2^{2p} ,2^{({\textstyle{1 \over 2}} - {\textstyle{p \over 2}} + 2p)} ,0,2^{( - {\textstyle{1 \over 2}} + {\textstyle{p \over 2}})} ,2))\,, 
 \end{eqnarray}
\begin{eqnarray}\label{equ.kxn883g}
&&\frac{\sqrt 3}{2}\log \left( {\frac{{(1 + 2^{^{( - {\textstyle{1 \over 2}} - {\textstyle{p \over 2}})} } \sqrt 3  - 2^{^{( - {\textstyle{1 \over 2}} - {\textstyle{p \over 2}})} }  + 2^{ - p} )^2 }}{{1 - 2^{^{({\textstyle{1 \over 2}} - {\textstyle{p \over 2}})} }  - 2^{^{({\textstyle{1 \over 2}} - {\textstyle{{3p} \over 2}})} }  + 2^{ - 2p}  + 2^{ - p} }}} \right)\nonumber\\
&&\quad= \frac{{3 }}{{2^{12p} }}P(1,2^{12p} ,24,(2^{( - {\textstyle{1 \over 2}} + {\textstyle{p \over 2}} + 11p)} ,2^{11p} ,0,0, \nonumber\\ 
&& \qquad 2^{( - {\textstyle{1 \over 2}} + {\textstyle{p \over 2}} + 9p)} ,0, - 2^{( - {\textstyle{1 \over 2}} + {\textstyle{p \over 2}} + 8p)} ,0,0, - 2^{7p} ,  \nonumber\\ 
&& \qquad - 2^{( - {\textstyle{1 \over 2}} + {\textstyle{p \over 2}} + 6p)} ,0, - 2^{( - {\textstyle{1 \over 2}} + {\textstyle{p \over 2}} + 5p)} , - 2^{5p} ,0,0,  \nonumber\\ 
&& \qquad - 2^{( - {\textstyle{1 \over 2}} + {\textstyle{p \over 2}} + 3p)} ,0,2^{( - {\textstyle{1 \over 2}} + {\textstyle{p \over 2}} + 2p)} ,0,0,2^p ,\nonumber\\
&& \qquad \qquad 2^{( - {\textstyle{1 \over 2}} + {\textstyle{p \over 2}})} ,0))\,, 
 \end{eqnarray}

\begin{eqnarray}\label{equ.dopqln2}
&&\sqrt 3\arctan \left[ {\left( {\frac{{1 - 2^{{\textstyle{{p + 1} \over 2}}} }}{{1 + 2^{{\textstyle{{p + 1} \over 2}}}  - 2^{p + 1} }}} \right)\sqrt 3 } \right]\nonumber\\
&&\quad= \frac{{3 }}{{2^{12p} }}P(1,2^{12p} ,24,(2^{( - {\textstyle{1 \over 2}} + {\textstyle{p \over 2}} + 11p)} , \nonumber\\ 
&&\qquad 0,0,2^{10p} ,2^{( - {\textstyle{1 \over 2}} + {\textstyle{p \over 2}} + 9p)} ,0,2^{( - {\textstyle{1 \over 2}} + {\textstyle{p \over 2}} + 8p)} ,2^{8p} ,0,0,2^{( - {\textstyle{1 \over 2}} + {\textstyle{p \over 2}} + 6p)} , \nonumber\\ 
&&\quad\qquad 0, - 2^{( - {\textstyle{1 \over 2}} + {\textstyle{p \over 2}} + 5p)} ,0,0, - 2^{4p} , - 2^{( - {\textstyle{1 \over 2}} + {\textstyle{p \over 2}} + 3p)} ,0, \nonumber\\ 
&&\qquad\qquad  - 2^{( - {\textstyle{1 \over 2}} + {\textstyle{p \over 2}} + 2p)} , - 2^{2p} ,0,0, - 2^{( - {\textstyle{1 \over 2}} + {\textstyle{p \over 2}})} ,0))\,,
 \end{eqnarray}

\begin{eqnarray}\label{equ.sbqy01x}
 - \arctan \left[ {\frac{{2^{{\textstyle{{1 - p} \over 2}}}  - 1}}{{ - 2^{{\textstyle{{1 + p} \over 2}}}  + 1}}} \right] &=& \frac{1}{{2^{12p} }}P(1,2^{12p} ,24,( - 2^{( - {\textstyle{1 \over 2}} + {\textstyle{p \over 2}} + 11p)} ,2^{11p} ,2^{({\textstyle{1 \over 2}} - {\textstyle{p \over 2}} + 11p)} , \nonumber\\ 
&& 0,2^{( - {\textstyle{1 \over 2}} + {\textstyle{p \over 2}} + 9p)} ,2^{1 + 9p} ,2^{( - {\textstyle{1 \over 2}} + {\textstyle{p \over 2}} + 8p)} ,0,2^{({\textstyle{1 \over 2}} - {\textstyle{p \over 2}} + 8p)} ,2^{7p} , - 2^{( - {\textstyle{1 \over 2}} + {\textstyle{p \over 2}} + 6p)} , \nonumber\\ 
&& 0,2^{( - {\textstyle{1 \over 2}} + {\textstyle{p \over 2}} + 5p)} , - 2^{5p} , - 2^{({\textstyle{1 \over 2}} - {\textstyle{p \over 2}} + 5p)} ,0, - 2^{( - {\textstyle{1 \over 2}} + {\textstyle{p \over 2}} + 3p)} , - 2^{1 + 3p} , \nonumber\\ 
&&  - 2^{( - {\textstyle{1 \over 2}} + {\textstyle{p \over 2}} + 2p)} ,0, - 2^{({\textstyle{1 \over 2}} - {\textstyle{p \over 2}} + 2p)} , - 2^p ,2^{( - {\textstyle{1 \over 2}} + {\textstyle{p \over 2}})} ,0))\,, 
 \end{eqnarray}

\begin{eqnarray}\label{equ.n2w2c23}
&& - \frac{1}{2}\log( {1 + 2^{^{({\textstyle{1 \over 2}} - {\textstyle{p \over 2}})} }  + 2^{ - p}  + 2^{^{({\textstyle{1 \over 2}} + {\textstyle{p \over 2}} - 2p)} }  + 2^{ - 2p} } ) \nonumber\\
&&= \frac{1}{{2^{12p} }}P\left( {1,2^{12p} ,24,\left( { - 2^{( - {\textstyle{1 \over 2}} + {\textstyle{p \over 2}} + 11p)} ,0,} \right.} \right. - 2^{({\textstyle{1 \over 2}} - {\textstyle{p \over 2}} + 11p)} ,2^{10p} , \nonumber\\ 
&& 2^{( - {\textstyle{1 \over 2}} + {\textstyle{p \over 2}} + 9p)} ,0, - 2^{( - {\textstyle{1 \over 2}} + {\textstyle{p \over 2}} + 8p)} , - 2^{8p} ,2^{({\textstyle{1 \over 2}} - {\textstyle{p \over 2}} + 8p)} ,0,2^{( - {\textstyle{1 \over 2}} + {\textstyle{p \over 2}} + 6p)} , \nonumber\\ 
&&  - 2^{1 + 6p} ,2^{( - {\textstyle{1 \over 2}} + {\textstyle{p \over 2}} + 5p)} ,0,2^{({\textstyle{1 \over 2}} - {\textstyle{p \over 2}} + 5p)} , - 2^{4p} , - 2^{( - {\textstyle{1 \over 2}} + {\textstyle{p \over 2}} + 3p)} ,0, \nonumber\\ 
&& \left. {\left. {2^{( - {\textstyle{1 \over 2}} + {\textstyle{p \over 2}} + 2p)} ,2^{2p} , - 2^{({\textstyle{1 \over 2}} - {\textstyle{p \over 2}} + 2p)} ,0, - 2^{( - {\textstyle{1 \over 2}} + {\textstyle{p \over 2}})} ,2} \right)} \right)\,, 
 \end{eqnarray}

\begin{eqnarray}
\frac{\sqrt 3}{2}\log \left[ {\frac{{1 + 2^{ - {\textstyle{1 \over 2}} - {\textstyle{p \over 2}}} (1 + \sqrt 3 ) + 2^{ - p} }}{{1 + 2^{ - {\textstyle{1 \over 2}} - {\textstyle{p \over 2}}} (1 - \sqrt 3 ) + 2^{ - p} }}} \right]
&=& \frac{{3 }}{{2^{12p} }}P\left( {1,2^{12p} ,24,\left( {2^{( - {\textstyle{1 \over 2}} + {\textstyle{p \over 2}} + 11p)} , - 2^{11p} ,} \right.} \right.0,0, \nonumber\\ 
&& 2^{( - {\textstyle{1 \over 2}} + {\textstyle{p \over 2}} + 9p)} ,0, - 2^{( - {\textstyle{1 \over 2}} + {\textstyle{p \over 2}} + 8p)} ,0,0,2^{7p} , - 2^{( - {\textstyle{1 \over 2}} + {\textstyle{p \over 2}} + 6p)} , \nonumber\\ 
&& 0, - 2^{( - {\textstyle{1 \over 2}} + {\textstyle{p \over 2}} + 5p)} ,2^{5p} ,0,0, - 2^{( - {\textstyle{1 \over 2}} + {\textstyle{p \over 2}} + 3p)} ,0, \nonumber\\ 
&& \left. {\left. {2^{( - {\textstyle{1 \over 2}} + {\textstyle{p \over 2}} + 2p)} ,0,0, - 2^p ,2^{( - {\textstyle{1 \over 2}} + {\textstyle{p \over 2}})} ,0} \right)} \right)\,, 
\end{eqnarray}

\begin{eqnarray}
\sqrt 3\arctan \left[ {\left( {\frac{{1 + 2^{{\textstyle{{p + 1} \over 2}}} }}{{-1 + 2^{{\textstyle{{p + 1} \over 2}}}  + 2^{p + 1} }}} \right)\sqrt 3 } \right] &=& \frac{{3 }}{{2^{12p} }}P\left( {1,2^{12p} ,24,\left( {2^{( - {\textstyle{1 \over 2}} + {\textstyle{p \over 2}} + 11p)} ,0,} \right.} \right.0, - 2^{10p} , \nonumber\\ 
&& 2^{( - {\textstyle{1 \over 2}} + {\textstyle{p \over 2}} + 9p)} ,0,2^{( - {\textstyle{1 \over 2}} + {\textstyle{p \over 2}} + 8p)} , - 2^{8p} ,0,0,2^{( - {\textstyle{1 \over 2}} + {\textstyle{p \over 2}} + 6p)} , \nonumber\\ 
&& 0, - 2^{( - {\textstyle{1 \over 2}} + {\textstyle{p \over 2}} + 5p)} ,0,0,2^{4p} , - 2^{( - {\textstyle{1 \over 2}} + {\textstyle{p \over 2}} + 3p)} ,0, \nonumber\\ 
&& \left. {\left. { - 2^{( - {\textstyle{1 \over 2}} + {\textstyle{p \over 2}} + 2p)} ,2^{2p} ,0,0, - 2^{( - {\textstyle{1 \over 2}} + {\textstyle{p \over 2}})} ,0} \right)} \right)\,, 
\end{eqnarray}

and

\begin{eqnarray}\label{equ.d9l4f7s}
\arctan \left[ {\frac{{2^{{\textstyle{{1 - p} \over 2}}}  + 1}}{{ 2^{{\textstyle{{1 + p} \over 2}}}  + 1}}} \right] &=& \frac{1}{{2^{12p} }}P\left( {1,2^{12p} ,24,\left( {2^{( - {\textstyle{1 \over 2}} + {\textstyle{p \over 2}} + 11p)} ,2^{11p} ,} \right.} \right. - 2^{({\textstyle{1 \over 2}} - {\textstyle{p \over 2}} + 11p)} ,0, \nonumber\\ 
&& - 2^{( - {\textstyle{1 \over 2}} + {\textstyle{p \over 2}} + 9p)} ,2^{1 + 9p} , - 2^{( - {\textstyle{1 \over 2}} + {\textstyle{p \over 2}} + 8p)} ,0, - 2^{({\textstyle{1 \over 2}} - {\textstyle{p \over 2}} + 8p)} ,2^{7p} ,2^{( - {\textstyle{1 \over 2}} + {\textstyle{p \over 2}} + 6p)} , \nonumber\\ 
&& 0, - 2^{( - {\textstyle{1 \over 2}} + {\textstyle{p \over 2}} + 5p)} , - 2^{5p} ,2^{({\textstyle{1 \over 2}} - {\textstyle{p \over 2}} + 5p)} ,0,2^{( - {\textstyle{1 \over 2}} + {\textstyle{p \over 2}} + 3p)} , - 2^{1 + 3p} , \nonumber\\ 
&& \left. {\left. {2^{( - {\textstyle{1 \over 2}} + {\textstyle{p \over 2}} + 2p)} ,0,2^{({\textstyle{1 \over 2}} - {\textstyle{p \over 2}} + 2p)} , - 2^p , - 2^{( - {\textstyle{1 \over 2}} + {\textstyle{p \over 2}})} ,0} \right)} \right)\,. 
 \end{eqnarray}

Particular cases of these formulas will be discussed in section~\ref{sec.nr81n5a}.

\section{Degree $2$ Formulas}

The imaginary part of the dilogarithm function can be expressed in closed form as~\cite{lewin81}

\begin{equation}\label{equ.ixjhxky}
 {\rm Im}\, {{\rm Li}_2 \left[ {qe^{ix} } \right]}  = \omega \log q + \frac{1}{2}{\rm Cl}_2 (2\omega ) - \frac{1}{2}{\rm Cl}_2 (2\omega  + 2x) + \frac{1}{2}{\rm Cl}_2 (2x ) \,, 
\end{equation}
where 
\[
\omega  = \arctan \left( {\frac{{q\sin x}}{{1 - q\cos x}}} \right)\,.
\]

Using \eqref{equ.ixjhxky}, the results \eqref{equ.aaob80z} and \eqref{equ.escpvp9} can be written in degree $2$ as

\begin{eqnarray}\label{equ.zf2kdeb}
&&  - (\omega _1  + \omega _2 )p\log 2 + {\rm Cl}_2 (2\omega _1 ) + {\rm Cl}_2 (2\omega _2 ) \nonumber\\ 
&&  - {\rm Cl}_2 (2\omega _1  + {\pi  \mathord{\left/
 {\vphantom {\pi  6}} \right.
 \kern-\nulldelimiterspace} 6}) + {\rm Cl}_2 ( - 2\omega _2  + 5{\pi  \mathord{\left/
 {\vphantom {\pi  6}} \right.
 \kern-\nulldelimiterspace} 6}) + {1 \mathord{\left/
 {\vphantom {1 2}} \right.
 \kern-\nulldelimiterspace} 2}{\rm Cl}_2 ({\pi  \mathord{\left/
 {\vphantom {\pi  3}} \right.
 \kern-\nulldelimiterspace} 3}) \nonumber\\ 
&&  = 2\,{\rm Im\,Li}_2 \left[ {\frac{1}{{\sqrt 2 ^p }}\exp \left( {\frac{{i\pi }}{{12}}} \right)} \right] + 2\,{\rm Im\,Li}_2 \left[ {\frac{1}{{\sqrt 2 ^p }}\exp \left( {\frac{{7i\pi }}{{12}}} \right)} \right] \nonumber\\
&&  = \frac{{\sqrt 3 }}{{2^{12p - 1} }}P(2,2^{12p} ,24,(2^{( - {\textstyle{1 \over 2}} + {\textstyle{p \over 2}} + 11p)} ,0,0,2^{10p} , \nonumber\\ 
&& 2^{( - {\textstyle{1 \over 2}} + {\textstyle{p \over 2}} + 9p)} ,0,2^{( - {\textstyle{1 \over 2}} + {\textstyle{p \over 2}} + 8p)} ,2^{8p} ,0,0,2^{( - {\textstyle{1 \over 2}} + {\textstyle{p \over 2}} + 6p)} , \nonumber\\ 
&& 0, - 2^{( - {\textstyle{1 \over 2}} + {\textstyle{p \over 2}} + 5p)} ,0,0, - 2^{4p} , - 2^{( - {\textstyle{1 \over 2}} + {\textstyle{p \over 2}} + 3p)} ,0, \nonumber\\ 
&&  - 2^{( - {\textstyle{1 \over 2}} + {\textstyle{p \over 2}} + 2p)} , - 2^{2p} ,0,0, - 2^{( - {\textstyle{1 \over 2}} + {\textstyle{p \over 2}})} ,0))
 \end{eqnarray}

and

\begin{eqnarray}\label{equ.k9injae}
&& (\omega _2  - \omega _1 )p\log 2 + {\rm Cl}_2 (2\omega _1 ) - {\rm Cl}_2 (2\omega _2 ) \nonumber\\ 
&&  - {\rm Cl}_2 (2\omega _1  + {\pi  \mathord{\left/
 {\vphantom {\pi  6}} \right.
 \kern-\nulldelimiterspace} 6}) - {\rm Cl}_2 ( - 2\omega _2  + 5{\pi  \mathord{\left/
 {\vphantom {\pi  6}} \right.
 \kern-\nulldelimiterspace} 6}) + {{4G} \mathord{\left/
 {\vphantom {{4G} 3}} \right.
 \kern-\nulldelimiterspace} 3} \nonumber\\ 
&&  = 2\,{\rm Im\,Li}_2 \left[ {\frac{1}{{\sqrt 2 ^p }}\exp \left( {\frac{{i\pi }}{{12}}} \right)} \right] - 2\,{\rm Im\,Li}_2 \left[ {\frac{1}{{\sqrt 2 ^p }}\exp \left( {\frac{{7i\pi }}{{12}}} \right)} \right] \nonumber\\ 
&&  = \frac{1}{{2^{12p - 1} }}P(2,2^{12p} ,24,( - 2^{( - {\textstyle{1 \over 2}} + {\textstyle{p \over 2}} + 11p)} ,2^{11p} ,2^{({\textstyle{1 \over 2}} - {\textstyle{p \over 2}} + 11p)} ,0, \nonumber\\ 
&& 2^{( - {\textstyle{1 \over 2}} + {\textstyle{p \over 2}} + 9p)} ,2^{1 + 9p} ,2^{( - {\textstyle{1 \over 2}} + {\textstyle{p \over 2}} + 8p)} ,0,2^{({\textstyle{1 \over 2}} - {\textstyle{p \over 2}} + 8p)} ,2^{7p} , - 2^{( - {\textstyle{1 \over 2}} + {\textstyle{p \over 2}} + 6p)} , \nonumber\\ 
&& 0,2^{( - {\textstyle{1 \over 2}} + {\textstyle{p \over 2}} + 5p)} , - 2^{5p} , - 2^{({\textstyle{1 \over 2}} - {\textstyle{p \over 2}} + 5p)} ,0, - 2^{( - {\textstyle{1 \over 2}} + {\textstyle{p \over 2}} + 3p)} , - 2^{1 + 3p} , \nonumber\\ 
&&  - 2^{( - {\textstyle{1 \over 2}} + {\textstyle{p \over 2}} + 2p)} ,0, - 2^{({\textstyle{1 \over 2}} - {\textstyle{p \over 2}} + 2p)} , - 2^p ,2^{( - {\textstyle{1 \over 2}} + {\textstyle{p \over 2}})} ,0))\,,
 \end{eqnarray}

where $\omega _1$ and $\omega _2$ are given by

\[
\tan \omega _1  = \frac{{\sqrt 3  - 1}}{{\sqrt 2 ^{p + 3}  - \sqrt 3  - 1}}
\]

and

\[
\tan \omega _2  = \frac{{\sqrt 3  + 1}}{{\sqrt 2 ^{p + 3}  + \sqrt 3  - 1}}\,.
\]

\bigskip

Similarly, using \eqref{equ.ixjhxky}, the results \eqref{equ.digtrb1} and \eqref{equ.uklodgq} can be written in degree~$2$ as
\begin{eqnarray}\label{equ.zf2kdeb}
&&  - (\omega _3  + \omega _4 )p\log 2 + {\rm Cl}_2 (2\omega _3 ) + {\rm Cl}_2 (2\omega _4 ) \nonumber\\ 
&&  - {\rm Cl}_2 (2\omega _3  + {5\pi  \mathord{\left/
 {\vphantom {\pi  6}} \right.
 \kern-\nulldelimiterspace} 6}) + {\rm Cl}_2 ( - 2\omega _4  + {\pi  \mathord{\left/
 {\vphantom {\pi  6}} \right.
 \kern-\nulldelimiterspace} 6}) - {1 \mathord{\left/
 {\vphantom {1 2}} \right.
 \kern-\nulldelimiterspace} 2}{\rm Cl}_2 ({\pi  \mathord{\left/
 {\vphantom {\pi  3}} \right.
 \kern-\nulldelimiterspace} 3}) \nonumber\\ 
&&  = 2\,{\rm Im\,Li}_2 \left[ {\frac{1}{{\sqrt 2 ^p }}\exp \left( {\frac{{5i\pi }}{{12}}} \right)} \right] + 2\,{\rm Im\,Li}_2 \left[ {\frac{1}{{\sqrt 2 ^p }}\exp \left( {\frac{{11i\pi }}{{12}}} \right)} \right] \nonumber\\
&& = \frac{{\sqrt 3 }}{{2^{12p-1} }}P\left( {2,2^{12p} ,24,\left( {2^{( - {\textstyle{1 \over 2}} + {\textstyle{p \over 2}} + 11p)} ,0,} \right.} \right.0, - 2^{10p} , \nonumber\\ 
&& 2^{( - {\textstyle{1 \over 2}} + {\textstyle{p \over 2}} + 9p)} ,0,2^{( - {\textstyle{1 \over 2}} + {\textstyle{p \over 2}} + 8p)} , - 2^{8p} ,0,0,2^{( - {\textstyle{1 \over 2}} + {\textstyle{p \over 2}} + 6p)} , \nonumber\\ 
&& 0, - 2^{( - {\textstyle{1 \over 2}} + {\textstyle{p \over 2}} + 5p)} ,0,0,2^{4p} , - 2^{( - {\textstyle{1 \over 2}} + {\textstyle{p \over 2}} + 3p)} ,0, \nonumber\\ 
&& \left. {\left. { - 2^{( - {\textstyle{1 \over 2}} + {\textstyle{p \over 2}} + 2p)} ,2^{2p} ,0,0, - 2^{( - {\textstyle{1 \over 2}} + {\textstyle{p \over 2}})} ,0} \right)} \right)
 \end{eqnarray}

and
\begin{eqnarray}\label{equ.k9injae}
&& (\omega _4  - \omega _3 )p\log 2 + {\rm Cl}_2 (2\omega _3 ) - {\rm Cl}_2 (2\omega _4 ) \nonumber\\ 
&&  - {\rm Cl}_2 (2\omega _3  + {5\pi  \mathord{\left/
 {\vphantom {\pi  6}} \right.
 \kern-\nulldelimiterspace} 6}) - {\rm Cl}_2 ( - 2\omega _4  + {\pi  \mathord{\left/
 {\vphantom {\pi  6}} \right.
 \kern-\nulldelimiterspace} 6}) + {{4G} \mathord{\left/
 {\vphantom {{4G} 3}} \right.
 \kern-\nulldelimiterspace} 3} \nonumber\\ 
&&  = 2\,{\rm Im\,Li}_2 \left[ {\frac{1}{{\sqrt 2 ^p }}\exp \left( {\frac{{5i\pi }}{{12}}} \right)} \right] - 2\,{\rm Im\,Li}_2 \left[ {\frac{1}{{\sqrt 2 ^p }}\exp \left( {\frac{{11i\pi }}{{12}}} \right)} \right] \nonumber\\ 
&&= \frac{1}{{2^{12p-1} }}P\left( {2,2^{12p} ,24,\left( {2^{( - {\textstyle{1 \over 2}} + {\textstyle{p \over 2}} + 11p)} ,2^{11p} ,} \right.} \right. - 2^{({\textstyle{1 \over 2}} - {\textstyle{p \over 2}} + 11p)} ,0, \nonumber\\ 
&& - 2^{( - {\textstyle{1 \over 2}} + {\textstyle{p \over 2}} + 9p)} ,2^{1 + 9p} , - 2^{( - {\textstyle{1 \over 2}} + {\textstyle{p \over 2}} + 8p)} ,0, - 2^{({\textstyle{1 \over 2}} - {\textstyle{p \over 2}} + 8p)} ,2^{7p} ,2^{( - {\textstyle{1 \over 2}} + {\textstyle{p \over 2}} + 6p)} , \nonumber\\ 
&& 0, - 2^{( - {\textstyle{1 \over 2}} + {\textstyle{p \over 2}} + 5p)} , - 2^{5p} ,2^{({\textstyle{1 \over 2}} - {\textstyle{p \over 2}} + 5p)} ,0,2^{( - {\textstyle{1 \over 2}} + {\textstyle{p \over 2}} + 3p)} , - 2^{1 + 3p} , \nonumber\\ 
&& \left. {\left. {2^{( - {\textstyle{1 \over 2}} + {\textstyle{p \over 2}} + 2p)} ,0,2^{({\textstyle{1 \over 2}} - {\textstyle{p \over 2}} + 2p)} , - 2^p , - 2^{( - {\textstyle{1 \over 2}} + {\textstyle{p \over 2}})} ,0} \right)} \right)\,,
 \end{eqnarray}

where $\omega _3$ and $\omega _4$ are given by
\[
\tan \omega _3 = \frac{{\sqrt 3  + 1}}{{\sqrt 2 ^{p + 3}  - \sqrt 3  + 1}}
\]

and
\[
\tan \omega _4  = \frac{{\sqrt 3  - 1}}{{\sqrt 2 ^{p + 3}  + \sqrt 3  + 1}}\,.
\]

In the above formulas, ${\rm G}={\rm Cl}_2(\pi/2)$ is Catalan's constant.

\bigskip

In deriving \eqref{equ.zf2kdeb} and \eqref{equ.k9injae}, we used \mbox{equation 4.32 pg 106} and \mbox{equation 4.17 pg 104} of \cite{lewin81}, namely,
\[
{\rm Cl}_2 \left( {\frac{\pi }{6}} \right) + {\rm Cl}_2 \left( {\frac{{5\pi }}{6}} \right) = \frac{{4{\rm G}}}{3}
\]

and
\[
{\rm Cl}_2 \left( x \right) - {\rm Cl}_2 \left( \pi-x \right) = \frac{1}{2}{\rm Cl}_2 \left( {2x} \right)\,.
\]

\section{Interesting Particular Cases}

\subsection{Degree $1$ Binary BBP-type Formulas}\label{sec.nr81n5a}
We first note that \eqref{equ.x9u9heu} and \eqref{equ.n2w2c23} give an infinite set of primes whose logarithms have binary BBP-type formulas, which serve to augment the known ones (e.g. those found in \cite{bailey09}) and \cite{chamberland03}. Similarly, \eqref{equ.sbqy01x} and \eqref{equ.d9l4f7s} give an infinite set of rationals whose arctangents have binary BBP-type representations. We now present a couple of binary degree $1$ formulas.

\subsubsection{Binary BBP-type formula for $\log 2$}

The identity
\begin{equation}\label{equ.fjvr030}
\log 2 = {\rm Re\,Li}_1 \left[ {\frac{1}{{\sqrt 2 }}\exp \left( {\frac{{i\pi }}{{12}}} \right)} \right] + {\rm Re\,Li}_1 \left[ {\frac{1}{{\sqrt 2 }}\exp \left( {\frac{{7i\pi }}{{12}}} \right)} \right]
\end{equation}

and $p=1$ in \eqref{equ.x9u9heu} lead to the binary BBP-type formula
\begin{eqnarray}
\log 2 &=&\frac{1}{2^{12}}P(1,2^{12},24,(2^{11}, 0, 2^{11}, 2^{10}, -2^9, 0, 2^8, -2^8,\nonumber\\
&&\quad-2^8, 0, -2^6, -2^7, -2^5, 0, -2^5, -2^4, 2^3, 0,\\
&&\qquad -2^2, 2^2, 2^2,0, 1, 2))\nonumber\,.
\end{eqnarray}

\subsubsection{Binary BBP-type formula for $\pi\sqrt 3$}

The identity

\begin{equation}\label{equ.vvt3cmq}
\frac{\pi }{3} = {\rm Im\,Li}_1 \left[ {\frac{1}{{\sqrt 2 }}\exp \left( {\frac{{i\pi }}{{12}}} \right)} \right] + {\rm Im\,Li}_1 \left[ {\frac{1}{{\sqrt 2 }}\exp \left( {\frac{{7i\pi }}{{12}}} \right)} \right]
\end{equation}

and $p=1$ in \eqref{equ.dopqln2} lead to the binary BBP-type formula

\begin{eqnarray}\label{equ.kd3ua8q}
\pi\sqrt 3 &=&\frac{9}{2^{12}}P(1,2^{12},24,(2^{11}, 0, 0, 2^{10}, 2^9, 0, 2^8, 2^8, 0, 0,\nonumber\\
&&\quad 2^6, 0, -2^5, 0, 0, -2^4, -2^3, 0, -2^2, -2^2, 0, 0, -1, 0))\,.
\end{eqnarray}

\subsubsection{Binary BBP-type formula for $\sqrt 3 \log (2 + \sqrt 3 )$}

The identity
\begin{equation}
\log (2 + \sqrt 3 ) = {\rm Re\,Li}_1 \left[ {\frac{1}{{\sqrt 2 }}\exp \left( {\frac{{i\pi }}{{12}}} \right)} \right] - {\rm Re\,Li}_1 \left[ {\frac{1}{{\sqrt 2 }}\exp \left( {\frac{{7i\pi }}{{12}}} \right)} \right]
\end{equation}
and $p=1$ in \eqref{equ.kxn883g} lead to the binary BBP-type formula
\begin{eqnarray}\label{equ.drs0vzb}
\sqrt 3 \,\log (2 + \sqrt 3 )&=&\frac{3}{2^{12}}P(1,2^{12},24,(2^{11}, 2^{11}, 0, 0, 2^9, 0, -2^8,\nonumber\\
&&\quad 0, 0, -2^7, -2^6, 0, -2^5, -2^5, 0, 0, -2^3,\nonumber\\
&&\qquad 0, 2^2, 0, 0, 2, 1, 0))
\end{eqnarray}

\subsubsection{Binary BBP-type formula for $\arctan 1/6$}

The identity

\begin{equation}
- \arctan \left( {\frac{1}{6}} \right) = {\rm Im\,Li}_1 \left[ {\frac{1}{{\sqrt 2 ^3 }}\exp \left( {\frac{{i\pi }}{{12}}} \right)} \right] - {\rm Im\,Li}_1 \left[ {\frac{1}{{\sqrt 2 ^3 }}\exp \left( {\frac{{7i\pi }}{{12}}} \right)} \right]
\end{equation}

and $p=3$ in \eqref{equ.sbqy01x} lead to the binary BBP-type formula
\begin{eqnarray}
\arctan \left( {\frac{1}{6}} \right)&=&\frac{1}{2^{35}}P(1,2^{36},24,(2^{33}, -2^{32}, -2^{31}, 0, -2^{27},\nonumber\\
&& -2^{27}, -2^{24}, 0, -2^{22}, -2^{20}, 2^{18}, 0, -2^{15}, 2^{14},\nonumber\\ 
&& 2^{13}, 0, 2^9, 2^9, 2^6, 0, 2^4, 2^2, -1, 0))\,.
\end{eqnarray}

\subsection{Degree $2$ Binary BBP-type Formulas}
When $p=1$ then $\omega _1=\omega _2=\pi/6$ and \eqref{equ.zf2kdeb} and \eqref{equ.k9injae} simplify to
\begin{equation}\label{equ.z2pmu0t}
- \frac{\pi }{3}\log 2 + \frac{5}{2}\,{\rm Cl}_2 \left( {\frac{\pi }{3}} \right) = 2\,{\rm Im\,Li}_2 \left[ {\frac{1}{{\sqrt 2 }}\exp \left( {\frac{{i\pi }}{{12}}} \right)} \right] + 2\,{\rm Im\,Li}_2 \left[ {\frac{1}{{\sqrt 2 }}\exp \left( {\frac{{7i\pi }}{{12}}} \right)} \right]
\end{equation}
and
\begin{equation}\label{equ.gps7q0j}
{\rm G} = 3\,{\rm Im\,Li}_2 \left[ {\frac{1}{{\sqrt 2 }}\exp \left( {\frac{{i7\pi }}{{12}}} \right)} \right] - 3\,{\rm Im\,Li}_2 \left[ {\frac{1}{{\sqrt 2 }}\exp \left( {\frac{{i\pi }}{{12}}} \right)} \right]\,,
\end{equation}
respectively.

\bigskip

Choosing $q=1/2$ and $x=\pi/3$ in \eqref{equ.ixjhxky} gives
\begin{equation}\label{equ.om15tau}
- \pi \log 2 + 5\,{\rm Cl}_2 \left( {\frac{\pi }{3}} \right) = 6\,{\rm Im\,Li}_2 \left[ {\frac{1}{2}\exp \left( {\frac{\pi }{3}} \right)} \right]\,.
\end{equation}

We note that
\begin{eqnarray}\label{equ.ngzuu4w}
 {\rm Im\,Li}_2 \left[ {\frac{1}{2}\exp \left( {\frac{{i\pi }}{3}} \right)} \right] &=& \sum\limits_{k = 1}^\infty  {\frac{1}{{2^k }}\frac{{\sin ({{k\pi } \mathord{\left/
 {\vphantom {{k\pi } 3}} \right.
 \kern-\nulldelimiterspace} 3})}}{{k^2 }}}  \nonumber\\ 
&&  =\frac{\sqrt 3}{2^{10}}\,P(2,2^{12},24,(0, 2^{10}, 0, 2^9, 0, 0, 0, -2^7, 0, -2^6, 0, 0, 0,\nonumber\\
&&\qquad 2^4, 0, 2^3, 0, 0, 0, -2, 0, -1, 0, 0))\,. 
 \end{eqnarray}

\subsubsection{Binary BBP-type formula for $\sqrt 3\,{\rm Cl}_2(\pi/3)$}
Eliminating $\pi\log 2$ between \eqref{equ.z2pmu0t} and \eqref{equ.om15tau} we have
\begin{eqnarray}\label{equ.cwe2yky}
 {\rm Cl}_2 \left( {\frac{\pi }{3}} \right) &=& \frac{{12}}{5}\,\left\{{\rm Im\,Li}_2 \left[ {\frac{1}{{\sqrt 2 }}\exp \left( {\frac{{i\pi }}{{12}}} \right)} \right] + {\rm Im\,Li}_2 \left[ {\frac{1}{{\sqrt 2 }}\exp \left( {\frac{{7i\pi }}{{12}}} \right)} \right]\right\} \nonumber\\ 
&&  - \frac{{12}}{5}\,{\rm Im\,Li}_2 \left[ {\frac{1}{2}\exp \left( {\frac{{i\pi }}{3}} \right)} \right]\,. 
 \end{eqnarray}
 
Using \eqref{equ.zf2kdeb} (with $p=1$) and \eqref{equ.ngzuu4w} in \eqref{equ.cwe2yky}, we obtain a binary BBP-type formula for $\sqrt 3\,{\rm Cl}_2(\pi/3)$:
\begin{eqnarray}\label{equ.hrsg8f7}
\sqrt 3 \,{\rm Cl}_2 \left( {\frac{\pi }{3}} \right) &=& \frac{9}{5\cdot 2^{10}}\,P(2,2^{12},24, (2^{11}, -2^{12}, 0, -2^{10}, 2^9,\nonumber\\
&&\quad 0, 2^8, 3\cdot 2^8, 0, 2^8, 2^6, 0, -2^5, -2^6,0,\nonumber\\
&&\quad -3\cdot2^4, -2^3, 0, -2^2, 2^2, 0, 2^2, -1, 0))\,.
\end{eqnarray}

\subsubsection{Binary BBP-type formula for $\pi\sqrt 3\log 2$}
Eliminating ${\rm Cl}_2(\pi/3)$ between \eqref{equ.z2pmu0t} and \eqref{equ.om15tau} we have
\begin{eqnarray}\label{equ.j4mh9xu}
 \pi\log 2 &=& 12\,\left\{{\rm Im\,Li}_2 \left[ {\frac{1}{{\sqrt 2 }}\exp \left( {\frac{{i\pi }}{{12}}} \right)} \right] + {\rm Im\,Li}_2 \left[ {\frac{1}{{\sqrt 2 }}\exp \left( {\frac{{7i\pi }}{{12}}} \right)} \right]\right\} \nonumber\\ 
&&  - 18\,{\rm Im\,Li}_2 \left[ {\frac{1}{2}\exp \left( {\frac{{i\pi }}{3}} \right)} \right]\,. 
 \end{eqnarray}

\bigskip
 
Using \eqref{equ.zf2kdeb} (with $p=1$) and \eqref{equ.ngzuu4w} in \eqref{equ.j4mh9xu}, we obtain a binary BBP-type formula for $\pi\sqrt 3\log 2$:
\begin{eqnarray}\label{equ.fa7mcws}
 \pi\sqrt 3\log 2 &=&\frac{9}{2^{10}} P(2,2^{12},24,(2^{11}, -3\cdot 2^{11}, 0, -2^{11}, 2^9, 0, 2^8, 2^{10}, 0, 3\cdot2^7, 2^6, 0,\nonumber\\
&& -2^5, -3\cdot2^5, 0, -2^6, -2^3, 0, -2^2, 2^3, 0, 6, -1, 0)\,.  
 \end{eqnarray}

A variant of formula \eqref{equ.fa7mcws} (formula~32, section~5, in the BBP-Compendium) was discovered experimentally over ten years ago, but is hitherto unproved. Formula~32 in the Compendium is now proved by adding rational multiples of two zero relations to \eqref{equ.fa7mcws} (See \cite{bailey09}).
 
\subsubsection{Binary BBP-type formula for Catalan's constant G}

\eqref{equ.gps7q0j} leads immediately to
\begin{eqnarray}\label{equ.oc6yqru}
 {\rm G} &=& \frac{3}{2^{12}}\,P(2,2^{12},24,(2^{11}, -2^{11}, -2^{11}, 0, -2^9, -2^{10},-2^8, 0,\nonumber\\
&& \qquad -2^8, -2^7, 2^6, 0, -2^5, 2^5, 2^5, 0, 2^3, 2^4, 2^2, 0, 2^2, 2, 1, 0))\,.
 \end{eqnarray}

\subsection{Binary Zero Relations}

The identity
\begin{equation}\label{equ.crtvn46}
\frac{\pi }{6} = {\rm Im\,Li}_1 \left[ {\frac{1}{2}\exp \left( {\frac{{i\pi }}{3}} \right)} \right] = \sum\limits_{k = 1}^\infty  {\frac{1}{{2^k }}\frac{{\sin ({{k\pi } \mathord{\left/
 {\vphantom {{k\pi } 3}} \right.
 \kern-\nulldelimiterspace} 3})}}{k}} 
\end{equation}
leads to the binary BBP-type formula
\begin{eqnarray}\label{equ.squ4uf4}
\pi\sqrt 3&=&\frac{9}{2^{10}}P(1,2^{12},24,(0, 2^{10}, 0, 2^9, 0, 0, 0, -2^7, 0,-2^6, 0, 0, 0, 2^4, 0,  \nonumber\\ 
&& \qquad 2^3, 0, 0, 0, -2, 0, -1, 0, 0))\,.
\end{eqnarray}

\bigskip

We note that this is the base $2^{12}$ version of the formula listed for $\pi\sqrt 3$ in section~4 of the BBP Compendium~\cite{bailey09}.

\bigskip

Combining \eqref{equ.vvt3cmq} and \eqref{equ.crtvn46}, we have the identity
\begin{eqnarray}
 0 &=& {\rm Im\,Li}_1 \left[ {\frac{1}{{\sqrt 2 }}\exp \left( {\frac{{i\pi }}{{12}}} \right)} \right] + {\rm Im\,Li}_1 \left[ {\frac{1}{{\sqrt 2 }}\exp \left( {\frac{{7i\pi }}{{12}}} \right)} \right] \nonumber\\ 
&&  - 2\,{\rm Im\,Li}_1 \left[ {\frac{1}{2}\exp \left( {\frac{{i\pi }}{3}} \right)} \right]\,, 
\end{eqnarray}

which, with the use of \eqref{equ.kd3ua8q} and \eqref{equ.squ4uf4} gives the binary zero relation
\begin{eqnarray}
0&=&P(1,2^{12},24,(2^{11}, -2^{12}, 0, -2^{10}, 2^9, 0, 2^8, 3\cdot 2^8, 0, 2^8, 2^6, 0, -2^5,\nonumber\\
&& -2^6, 0, -3\cdot 2^4, -2^3, 0, -2^2, 2^2, 0, 2^2, -1, 0))\,.
\end{eqnarray}

\eqref{equ.fjvr030} and the formula ${\rm Li}_1 (1/2) = \log 2$ give the identity
\begin{equation}
0 = {\rm Re\,Li}_1 \left[ {\frac{1}{{\sqrt 2 }}\exp \left( {\frac{{i\pi }}{{12}}} \right)} \right] + {\rm Re\,Li}_1 \left[ {\frac{1}{{\sqrt 2 }}\exp \left( {\frac{{7i\pi }}{{12}}} \right)} \right] - {\rm Li}_1 \left[ {\frac{1}{2}} \right]
\end{equation}
which leads to the binary BBP-type zero relation
\begin{eqnarray}
0&=&P(1,2^{12},24,(2^{11}, -2^{12}, 2^{11}, -2^{10}, -2^9, -2^{10}, 2^8, -3\cdot 2^8, -2^8, -2^8, -2^6,\nonumber\\
&& -2^8, -2^5, -2^6, -2^5, -3\cdot 2^4, 2^3, -2^4, -2^2, -2^2, 2^2, -2^2, 1, 0))\,.
\end{eqnarray}

The identity
\begin{equation}
0 = {\rm Im\,Li}_1 \left[ {\frac{1}{{\sqrt 2 }}\exp \left( {\frac{{i\pi }}{{12}}} \right)} \right] - {\rm Im\,Li}_1 \left[ {\frac{1}{{\sqrt 2 }}\exp \left( {\frac{{7i\pi }}{{12}}} \right)} \right]
\end{equation}

and $p=1$ in \eqref{equ.sbqy01x} lead to the binary BBP-type zero relation
\begin{eqnarray}\label{equ.csgyn9z}
0&=&P(1,2^{12},24,(-2^{11}, 2^{11}, 2^{11}, 0, 2^9, 2^{10}, 2^8, 0, 2^8, 2^7, -2^6, 0, 2^5, -2^5,\nonumber\\ 
&& -2^5, 0, -2^3, -2^4, -2^2, 0, -2^2, -2, 1, 0))\,.
\end{eqnarray}

\bigskip

The identity
\begin{equation}
0 = {\rm Re\,Li}_1 \left[ {\frac{1}{{\sqrt 2 }}\exp \left( {\frac{{i\pi }}{{12}}} \right)} \right] + {\rm Re\,Li}_1 \left[ {\frac{1}{{\sqrt 2 }}\exp \left( {\frac{{7i\pi }}{{12}}} \right)} \right] - 2\,{\rm Re\, Li}_1 \left[ {\frac{1}{{\sqrt 2 }}\exp \left( {\frac{{i\pi }}{4}} \right)} \right]
\end{equation}

leads to the binary zero relation
\begin{eqnarray}
0&=&P(1,2^{12},24,(2^{11}, 0, -2^{12}, -3\cdot 2^{10}, -2^9, 0, 2^8, 3\cdot 2^8, 2^9, 0, -2^6, 0,\nonumber\\
&& -2^5, 0, 2^6,3\cdot 2^4, 2^3, 0, -2^2, -3\cdot 2^2, -2^3, 0, 1, 0))\,.
\end{eqnarray}

The identity
\begin{eqnarray}
 0 &=& {\rm Re\,Li}_1 \left[ {\frac{1}{{\sqrt 2 }}\exp \left( {\frac{{i\pi }}{{12}}} \right)} \right] + {\rm Re\,Li}_1 \left[ {\frac{1}{{\sqrt 2 }}\exp \left( {\frac{{7i\pi }}{{12}}} \right)} \right] \nonumber\\ 
&&  + 2\,{\rm Re\,Li}_1 \left[ {\frac{1}{{\sqrt 2 }}\exp \left( {\frac{{i\pi }}{2}} \right)} \right] - 2\,{\rm Re\,Li}_1 \left[ {\frac{1}{2}\exp \left( {\frac{{i\pi }}{3}} \right)} \right] 
 \end{eqnarray}

leads to the binary zero relation
\begin{eqnarray}
0&=&P(1,2^{12},24,(2^{11}, -2^{13}, 2^{11}, 5\cdot 2^{10}, -2^9, 2^{10}, 2^8, 3\cdot 2^8, -2^8, -2^9, -2^6, -2^8, -2^5,\nonumber\\
&& -2^7, -2^5,3\cdot 2^4, 2^3, 2^4, -2^2, 5\cdot 2^2, 2^2, -2^3, 1, 0))\,.
\end{eqnarray}

More degree~1 binary BBP-type zero relations in base $2^{12}$ and other bases can also be found in~\cite{lafont2011}.

\section{Conclusion}
Using a clear and straightforward approach, explicit digit extraction BBP-type formulas in very general binary bases were discovered. As particular examples, new binary formulas were obtained for $\pi\sqrt 3$, $\sqrt 3\pi\log 2$ and some other polylogarithm constants. New binary BBP-type zero relations were also established.

\end{document}